\documentclass[10pt,a4paper]{article}  

\usepackage{scalerel} 

\usepackage{parskip}

\usepackage[top=2.0cm,bottom=1.9cm,left=1.5cm,right=1.5cm]{geometry} 

\usepackage{multicol} 

\usepackage[authoryear,round,elide,nonamebreak]{natbib}

\usepackage[utf8]{inputenc}
\usepackage[french,english]{babel}

\usepackage[T1]{fontenc}

\usepackage{lmodern}

\usepackage[normalem]{ulem}

\newcommand\keyw{\textsl} 

\usepackage{paralist} 
\setdefaultitem{$\circ$}{$\triangleright$}{$\star$}{-} 
\usepackage{enumitem} 
\setlist{topsep=2pt,itemsep=2pt,partopsep=2pt,parsep=2pt,leftmargin=1em}

\usepackage[svgnames]{xcolor} 

\usepackage{amsfonts,amssymb}    
\usepackage[intlimits]{amsmath}  
\usepackage[]{mathalfa}
\usepackage{mathrsfs}            
\usepackage{upgreek}             
\usepackage{dsfont}              
\usepackage{mathtools}           
\usepackage{relsize}             
\usepackage{cancel}              
\usepackage{accents}             
\usepackage{bm}                  

\usepackage{float}      
\usepackage[small,bf]{caption} 
\usepackage{subfig} 

\usepackage{mdframed} 
\usepackage{csquotes} 

\usepackage{booktabs}   

\usepackage{amsthm}                                                  

\usepackage{algorithm}
\usepackage{algorithmic}

\definecolor{mycc}{rgb}{0, 0, 0.45}   
\definecolor{mygreen}{rgb}{0, 0.4, 0} 

\usepackage{ifpdf}                    
\ifpdf                                
	\usepackage[pdftex,               
		pdftitle={},                  
		pdfauthor={},                 
		hyperindex=true,              
		colorlinks=true,              
		linkcolor=mygreen,            
		citecolor=mycc,               
        urlcolor=blue,                
        ]{hyperref}                   
	\usepackage{epstopdf}             

	\usepackage[stretch=10]{microtype} 
\else                                  
	\usepackage[hypertex,              
		hyperindex=true,               
		colorlinks=false]{hyperref}    
	\usepackage{graphicx}              
\fi                                    
\usepackage{hyperref}

\usepackage[doipre={DOI:~}]{uri} 

\graphicspath{{./}{Pics/}}

\usepackage{cleveref}
\crefname{equation}{equation}{equations}
\Crefname{equation}{Equation}{Equations}
\crefname{subeqns}{equation}{equations}
\Crefname{subeqns}{Equation}{Equations}
\crefformat{subeqns}{equation~(#2#1#3)}
\Crefformat{subeqns}{Equation~(#2#1#3)}
\labelcrefformat{subeqns}{(#2#1#3)}
\crefname{figure}{Figure}{Figures}
\Crefname{figure}{Figure}{Figures}
\crefformat{subfigure}{Figure~(#2#1#3)}
\Crefformat{subfigure}{Figure~(#2#1#3)}
\crefname{subfigure}{Figure}{Figures}
\Crefname{subfigure}{Figure}{Figures}
\crefname{table}{Table}{Tables}
\Crefname{table}{Table}{Tables}
\crefname{theo}{Theorem}{Theorems}
\Crefname{theo}{Theorem}{Theorems}
\crefname{theo}{Theorem}{Theorems}
\Crefname{theo}{Theorem}{Theorems}
\crefname{lemm}{Lemma}{Lemmas}
\Crefname{lemm}{Lemma}{Lemmas}
\crefname{conj}{Conjecture}{Conjectures}
\Crefname{conj}{Conjecture}{Conjectures}
\crefname{defi}{Definition}{Definitions}
\Crefname{defi}{Definition}{Definitions}
\crefname{prop}{Proposition}{Propositions}
\Crefname{prop}{Proposition}{Propositions}
\crefname{property}{Property}{Properties}
\Crefname{property}{Property}{Properties}
\crefname{coro}{Corollary}{Corollaries}
\Crefname{coro}{Corollary}{Corollaries}
\crefname{algocf}{Algorithm}{Algorithms}
\Crefname{algocf}{Algorithm}{Algorithms}

\newcommand{\mat}[1]{{\mathbf{{#1}}}}
\newcommand{\vect}[1]{{\bm{#1}}}

\newcommand{\ones}[0]{\mathds{1}}

\newcommand{\tn}[1]{{\textnormal{{#1}}}} 

\newcommand{\Reals}{\mathbb{R}}

\DeclareMathOperator{\Expect}{\mathbb{E}}
\DeclareMathOperator{\Var}{\mathbb{V}\tn{ar}}
\DeclareMathOperator{\Cov}{\mathbb{C}\tn{ov}}

\newcommand{\NormDist}{\mathcal{N}}

\newcommand{\TP}{{\mathsf{T}}}
\newcommand{\tr}{\ensuremath{^{\TP}}}
\newcommand{\pinvsign}{{+}}
\newcommand{\pinv}{\ensuremath{^\pinvsign}}

\newcommand{\pd}[2]{\frac{\mathrm{\partial}{#1}}{\mathrm{\partial}{#2}}}

\newcommand{\Loss}[0]{\ensuremath{\mathcal{L}}}
\newcommand{\loss}[0]{\ensuremath{\ell}}
\newcommand{\ensgrad}[0]{\ensuremath{\bar{\nabla}}}
\newcommand{\Ut}[0]{\ensuremath{\tilde{\mat{U}}}}

\defcitealias{fonseca2017stochastic}{Fons17}
\defcitealias{chen2009efficient}{Chen09}
\newcommand{\Fonseca}[0]{\citetalias{fonseca2017stochastic}}
\newcommand{\Chen}[0]{\citetalias{chen2009efficient}}

\begin{document}

\newcommand{\todo}[1]{{{{\color{Maroon}{\textbf{TODO:} #1}}}}}

\title{
    \vspace{-1cm}
    {
        Review of ensemble gradients for robust optimisation
    }
}
\author{Patrick N. Raanes, Andreas S. Stordal, Rolf J. Lorentzen}
\date{\today}
\maketitle

\begin{abstract}
    In robust optimisation problems the objective function
    consists of an average over (an ensemble of) uncertain parameters.
    Ensemble optimisation (EnOpt) implements steepest descent
    by estimating the gradient using linear regression on Monte-Carlo simulations
    of (an ensemble of) control parameters.
    Applying EnOpt for robust optimisation is costly
    unless the evaluations over the two ensembles are combined, i.e. \enquote*{paired}.
    Here, we provide a new and more rigorous perspective on the
    stochastic simplex approximate gradient (StoSAG) used in EnOpt,
    explaining how it addresses detrimental cross-correlations arising from pairing
    by only capturing the variability due to the control vector,
    and not the vector of uncertain parameters.
    A few minor variants are derived from a generalized derivation,
    as well as a new approach using decorrelation.
    These variants are tested on linear and non-linear toy gradient estimation problems,
    where they achieve highly similar accuracy,
    but require a very large ensemble size to outperform the non-robust approach
    when accounting for variance and not just bias.
    Other original contributions include
    a discussion of the particular robust control objectives for which EnOpt is suited,
    illustrations,
    a variance reduction perspective,
    and a discussion on the centring in covariance and gradient estimation.
\end{abstract}

\begin{multicols}{2}

\section{Introduction}%
\label{sec:intro}

Ensemble methods are particularly skilled for estimation and control problems
involving dynamics that are nonlinear, high-dimensional, and computationally demanding,
e.g. fluid dynamics.
In data assimilation, the most prominent ensemble method
is the \keyw{ensemble Kalman filter} \citep[EnKF,][]{evensen1994sequential},
or some iterative version of it \citep[][]{bocquet2014iterative}.
The EnKF resembles the \keyw{unscented Kalman filter} \citep[][]{julier2004unscented},
but places a larger emphasis on using a reduced number of members, or particles,
and is usually complemented with techniques such as \keyw{localisation}
to compensate for the limited sample size.
In history matching and inverse problems it is the \keyw{iterative ensemble smoother}
\citep[EnRML or IES,][]{oliver1996conditional,raanes2019revising}.
In robust control problems it is
\keyw{ensemble optimisation} \citep[EnOpt,][hereafter \Chen]{chen2009efficient},
our focus here.
Common to all of these ensemble methods is that they
combine (non-invasive and embarrassingly parallelizable)
Monte-Carlo simulations with linear least-squares (LLS) regression to estimate
derivatives 
which are generally not readily available.

\subsection{Iterative convergence in the presence of stochasticity}%
\label{sub:CV}
Being that it uses a stochastic gradient for steepest descent/ascent,
EnOpt may be said fall in the category of \keyw{stochastic approximation}
\citep[SA,][]{robbins1951stochastic},
whose main result is a stochastic form of the contraction theorem:
even if we only have access to noisy, i.e. stochastic, function evaluations,
fixed-point iterations may converge,
providing the noise is unbiased and bounded,
and that the step length goes to zero at a specific rate.
\citet{kiefer1952stochastic} applied this result for optimisation
by considering the first-order, central finite difference approximations
for the gradient of an objective which again can only be evaluated subject to noise,
The cost of the evaluations required by this
\keyw{finite-difference stochastic approximation} (FDSA)
may be offset by parallelisation, but can be significant.

By contrast,
the simultaneous perturbation stochastic approximation method
\citet[SPSA][]{spall1992multivariate} 
uses only a single perturbation for its finite difference approximation,
only checking which of the two opposite directions are uphill or downhill,
similarly to the method of coordinate descent.
SPSA generally converges at a slower rate than FDSA, but at a lower cost.
Since the perturbation direction changes with each iteration,
some convergence results can be proven \citep[][\S 5.2.3]{bhatnagar2013stochastic}.

EnOpt may be said to fall somewhere in between FDSA and SPSA,
since we are generally free to choose the size
of the ensemble of controls to which LLS regression is applied.
Thus, despite the gradient estimates being stochastic
and not exploring all directions at each iteration,
we have some confidence that convergence
may be achieved in fairly general circumstances.
The remainder of this paper is not concerned iterative convergence,
focusing instead on the accuracy of the gradient estimates for a given ensemble.

\subsection{Outline}%
\label{sub:outline}
\Cref{sec:ens} begins by a discussion of covariance estimation
and centring strategies, which is complemented by \cref{sec:centring}.
We then define the fundamental regression estimator,
linking regression to average, analytic derivatives,
illustrating that it effectively blurs the objective (or cost) function,
and discussing common preconditioning techniques.
\Cref{sec:robust_control} continues the introduction,
detailing the type of problem at which EnOpt excels,
namely robust control problems where the objective function
involves computationally demanding simulations.
This specialised skill is largely due, as described in \cref{sec:reducingN},
to the cost-saving technique of \Chen\ called \enquote*{pairing},
and the refinement thereof of \citet{fonseca2014robust} called StoSAG,
which greatly improves its accuracy.

We think that a clear derivation of StoSAG is missing from the literature,
as discussed in \cref{sec:fonseca},
and that the surrounding theory can be clarified,
which is our main purpose and original contribution.
Our first explanation of the refinement is presented in \cref{sec:decoupling},
which is complemented by a variance reduction perspective in \cref{sec:variance_reduction},
and a linear example in \cref{sec:linear_example}.
\Cref{sec:general_stosag} derives some minor variations on the StoSAG gradient estimate,
all of which are experimentally compared in \cref{sec:experiments}.

\section{Fundamental ensemble gradient}%
\label{sec:ens}
An \keyw{ensemble} is a sample for which each
\keyw{realisation}, or \keyw{member},
$\vect{u}_n \in \Reals^{d_{\vect{u}}}$, for $n = 1, \ldots, N$,
is assumed an independent draw from the same probability distribution,
$p(\vect{u})$.

\subsection{Covariance estimation and centring}%
\label{sub:cov}

Denote $\bar{\vect{u}} = \frac{1}{N} \sum_n \vect{u}_n$ the sample mean,
where the summation bounds are left implicit.
Consider the cross-covariance of $\vect{u}$ with the generic function
$\vect{f}: \Reals^{d_{\vect{u}}} \rightarrow \Reals^{d_{\vect{f}}}$,
denoted $\mat{C}_{\vect{f}, \vect{u}} \in \Reals^{d_{\vect{f}} \times d_{\vect{u}}}$.
For $N \ge 2$, an unbiased estimate is
\begin{equation}
    \label{eq:ens_cov}
    \bar{\mat{C}}_{\vect{f}, \vect{u}}
    \coloneqq \frac{1}{N-1} \sum_n \vect{f}(\vect{u}_n)\, (\vect{u}_n - \bar{\vect{u}})\tr
    \,,
\end{equation}
where vectors have column orientation by default.
Note that, for any $c$ constant in $n$,
$\sum_n c \, (\vect{u}_n - \bar{\vect{u}}) = \vect{0}$,
so that the explicit \keyw{centring} of $\vect{f}(\vect{u}_n)$ in \cref{eq:ens_cov}
by subtracting its sample mean, $\frac{1}{N} \sum_n \vect{f}(\vect{u}_n)$,
can be omitted.

Unlike the data assimilation setting,
the true mean of the \keyw{input}, $\vect{u}$,
denoted $\vect{\mu}$, is known in EnOpt.
Therefore it is possible to estimate $\mat{C}_{\vect{f}, \vect{u}}$,
even for $N = 1$, by
\begin{equation}
    \label{eq:cov_estim_true_mean}
    \frac{1}{N} \sum_n \vect{f}(\vect{u}_n) (\vect{u}_n - \vect{\mu})\tr 
    \,.
\end{equation}
Even though $\vect{f}(\vect{u}_n)$ does not get subtracted
by the (unknown) true mean of the output,
$\Expect \vect{f}(\vect{u})$,
the estimate \eqref{eq:cov_estim_true_mean} is also unbiased,
since $\Expect b \, (\vect{u}_n - \vect{\mu}) = \vect{0}$
for any constant $b$.
Nevertheless, despite the seeming benefit of $\vect{\mu}$ over $\bar{\vect{u}}$,
the estimate of expression \eqref{eq:cov_estim_true_mean}
generally has a higher variance than $\bar{\mat{C}}_{\vect{f},\vect{u}}$ of \cref{eq:ens_cov}.
This is proven for linear $\vect{f}$ in \cref{sec:centring},
whose analysis also suggests
approximately centring $\vect{f}(\vect{u}_n)$ in expression \eqref{eq:cov_estim_true_mean}
by subtracting $\vect{f}(\vect{\mu})$,
even though $\vect{f}(\vect{\mu}) \neq \Expect \vect{f}(\vect{u})$ unless $\vect{f}$ is linear.

However, in our numerical tests (not shown),
the lowest variance is usually achieved
if the input ensemble, $\{\vect{u}_n\}$, is centred
\emph{before} running the Monte-Carlo simulation altogether,
such that $\bar{\vect{u}} = \vect{\mu}$.
Therefore, \emph{this operation will be tacitly assumed for the remainder of this manuscript.}

Let $\mat{U} \in \Reals^{d_{\vect{u}} \times N}$ be the ensemble matrix,
meaning that its $n$-th \emph{column} is the $n$-th member, $\vect{u}_n$,
and let $\Ut$ be the \emph{centred} ensemble matrix,
whose $n$-th column is the anomaly, or deviation, $\vect{u}_n - \bar{\vect{u}}$,
where, we recall, $\bar{\vect{u}} = \vect{\mu}$.
Also extend the definition of $\vect{f}$ to allow matrix input
by applying the original function to each column.
Then the sample covariance matrix \eqref{eq:ens_cov} may be written as
\begin{equation}
    \label{eq:ens_cov2}
    \bar{\mat{C}}_{\vect{f}, \vect{u}}
    = \frac{1}{N-1} \vect{f}(\mat{U})\, \Ut\tr
    \,.
\end{equation}

\subsection{Linear regression}%
\label{sub:linear_regression}
For a given distribution of $\vect{u}$,
the \keyw{(exact) regression coefficient matrix} of $\vect{f}$ is defined as
$\mat{C}_{\vect{f}, \vect{u}} \, \mat{C}_{\vect{u}}^{-1}$.
For any $N \ge 2$,
an estimate can be composed
using the covariance estimator \eqref{eq:ens_cov2}
and the (Moore-Penrose) pseudo-inverse
as 
$\ensgrad\! \vect{f} \coloneqq \bar{\mat{C}}_{\vect{f}, \vect{u}}\, \bar{\mat{C}}_{\vect{u}}\pinv$.
But, by the \emph{identity}
\begin{equation}
    \label{eq:pinv_identity}
    \mat{A}\tr (\mat{A} \mat{A}\tr)\pinv = \mat{A}\pinv
    \,,
\end{equation}
the regression estimate formula simplifies to
\begin{equation}
    \label{eq:ensgrad}
    \ensgrad\! \vect{f}
    =
    \vect{f}(\mat{U})\, \Ut\pinv
    \,.
\end{equation}
The right-hand side can be recognised as
the linear least squares (LLS) coefficient matrix,
meaning that it minimises the sum of squared residuals
among all linear approximations to $\vect{f}$.
Assuming Gaussianity of the residuals makes it a maximum-likelihood (ML) estimate.
Under other qualifying assumptions,
it can also be termed
BLUE, (U)MVUE, and minimum-norm LLS estimate.

\Cref{eq:ensgrad} is meant as a definition of the \emph{operator} $\ensgrad$.
Since $\ensgrad$ also entails Monte-Carlo simulations to generate $\vect{f}(\mat{U})$,
it may be useful to distinguish this notion from the more general LLS \emph{method},
which relates any two samples, not necessarily computed via some $\vect{f}$,
and we therefore call $\ensgrad\! \vect{f}$ the \keyw{ensemble (regression) gradient}
(or Jacobian). 
\citet{gu2007iterative} called $\ensgrad\! \vect{f}$
the \keyw{ensemble sensitivity matrix}.
\citet[][]{do2013theoretical} also discussed the closely-related \keyw{simplex gradient},
the difference being that the deviations are not centred,
but measured from a chosen reference point,
typically, but not necessarily, $\vect{\mu}$ and $\vect{f}(\vect{\mu})$.

\subsection{Asymptotics and Stein's lemma}%
\label{sub:stein}

Although both $\bar{\mat{C}}_{\vect{\vect{f}}, \vect{u}}$ and $\bar{\mat{C}}_{\vect{u}}$ are unbiased,
the ensemble regression gradient, $\ensgrad\! \vect{f}$ of \cref{eq:ensgrad},
does not generally have expectation $\mat{C}_{\vect{f}, \vect{u}}\, \mat{C}_{\vect{u}}^{-1}$,
unless $\vect{f}$ is linear and $N > d_{\vect{u}}$.
Nevertheless, by Slutsky's theorem, it \emph{is} consistent,
converging almost surely:
\begin{equation}
    \label{eq:ensgrad_consistent}
    \ensgrad\! \vect{f} \xrightarrow[N \rightarrow +\infty]{} \mat{C}_{\vect{f}, \vect{u}}\, \mat{C}_{\vect{u}}^{-1}
    \,.
\end{equation}

Now, suppose $\vect{u}$ is Gaussian,
$\vect{u} \sim \NormDist(\vect{\mu}, \mat{C}_{\vect{u}})$,
whose density is denoted $p(\vect{u}) \propto \exp(-\| \vect{u} - \vect{\mu} \|^2_{\mat{C}_u} / 2)$,
and that $\vect{f}$ is a function with bounded gradient, $\nabla\! \vect{f}(\vect{u})$.
Then
\begin{equation}
    \label{eq:cov_grad}
    \Expect \nabla\! \vect{f}(\vect{u})
    = \mat{C}_{\vect{f}, \vect{u}}\, \mat{C}_{\vect{u}}^{-1} \,,
\end{equation}
which is known as Stein's lemma \citep{stein1972bound}.
It may be shown by performing integration by parts for the expectation,
followed by Stokes' theorem and 
the fact that $p(\vect{u}) \rightarrow 0$ at infinity,
followed by the \enquote{logarithm trick},
i.e. $\nabla p(\vect{u}) = p(\vect{u}) \nabla\! \log p(\vect{u})$,
where Gaussianity yields
$\nabla\! \log p(\vect{u}) = - \mat{C}_{\vect{u}}^{-1} [\vect{u} - \vect{\mu}]$.

\Cref{eq:ensgrad_consistent,eq:cov_grad} produce
\begin{equation}
    \label{eq:ensgrad_lim}
    \ensgrad\! \vect{f}
    \xrightarrow[N]{}
    \Expect \nabla\! \vect{f}(\vect{u})
    \,,
\end{equation}
which establishes a relationship
from ensemble sensitivity to the actual gradient,
motivating its use in optimisation contexts or perspectives.
This is usually accomplished in the ensemble literature by inserting
the leading-order Taylor expansion into the cross-covariance.
However, the above approach by Stein's lemma should be favoured because it
(i) is exact,
(ii) shows that errors cancel thanks to averaging,
(iii) shows that the \enquote*{target} of the estimation is the
\emph{average gradient} -- not the gradient at a point.

\subsection{Implicit Gaussian blurring}%
\label{sub:as_kernel_smoothing}

\Cref{eq:ensgrad_lim} was initially pointed out in the ensemble literature by
\citet[][]{raanes2019revising}.
A similar result was found by \citet[][]{stordal2016theoretical} 
in terms of \keyw{mutation} gradient,
$\nabla_{\!\vect{\mu}} \Expect \vect{f}(\vect{u})$,
i.e. the gradient in a parameter of the distribution.
The two are related by the fact that, for any density
that can be written $p(\vect{u}) = p_0(\vect{u} - \vect{\theta})$,
meaning that $\vect{\theta}$ is a \keyw{location} parameter,
\begin{equation}
    \label{eq:exp_grad_exp}
    \nabla_{\!\vect{\theta}} \Expect \vect{f}(\vect{u})
    =
    \Expect \nabla\! \vect{f}(\vect{u})
    \,,
\end{equation}
whose demonstration involves integration by parts as in \cref{eq:cov_grad},
and requires that the \keyw{kernel}, $p_0$, be sufficiently decaying.

\begin{figure}[H]
    \centering
    \includegraphics[width=0.99\linewidth]{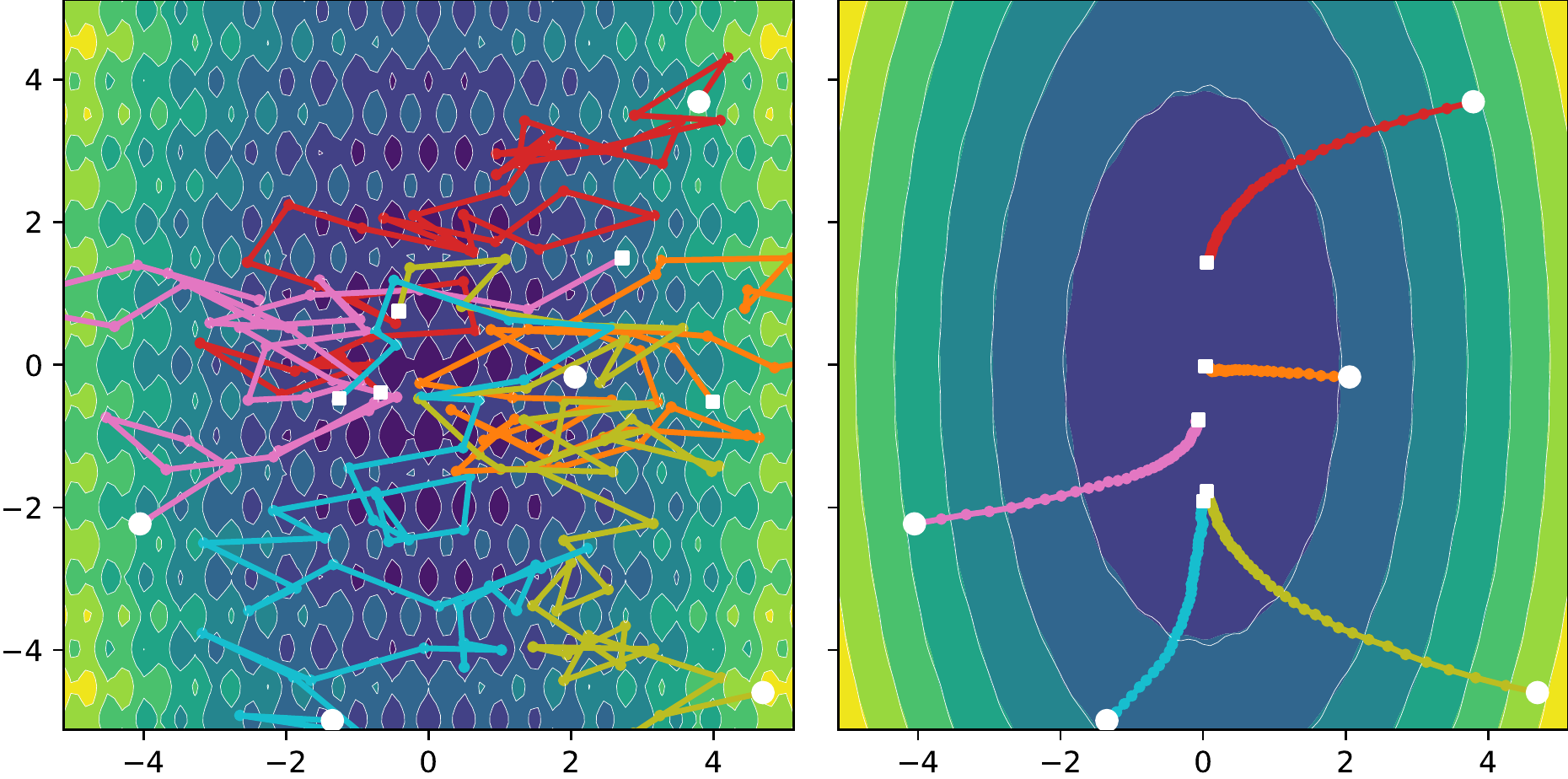}
    \caption{
        Steepest descent,
        using exact (analytic) gradient (left)
        and (asymptotic) ensemble gradient (right),
        illustrated on the Rastrigin objective,
        $\Loss(\vect{u}) = 20 + \sum_{i=1,2} \bigl[v_i^2 - 10 \cos(2 \pi v_i)\bigr]$,
        with $v_1 = 2 u_1$ introducing a degree of non-homogeneity
        and $v_2 = u_2$.
        The objective is shown on the left as contour surfaces ranging from 0 to 154,
        displaying its many local minima.
        The contours on the right show the blurred objective function,
        $\vect{\mu} \mapsto \Expect \Loss(\vect{u})$,
        with $\vect{u} \sim \NormDist(\vect{\mu}, \mat{I})$,
        which is what using the ensemble gradient implicitly optimises.
        Steps are taken as $0.012$ times the gradient (or its average),
        starting from $5$ different initial positions, shown as white disks;
        smaller steps cause the gradient descent with exact gradient to be stuck
        within horizontal layers.
        In practice, acceleration techniques (e.g. momentum) should be used,
        allowing for much smaller ensemble or larger system size.
    }
    \label{fig:rastrigin}
\end{figure}

Thus, if $\Loss: \Reals^{d_{\vect{u}}} \rightarrow \Reals$
is the objective to be minimised,
gradient descent using $\ensgrad \Loss$,
i.e. EnOpt, can be said to optimise, with respect to $\vect{\mu}$,
the smoothed, or \enquote*{blurred} objective function,
so-called because it is a convolution with the Gaussian kernel.
Optimisation of kernel-smoothed objectives were studied at least as early as
\citet{katkovnik1972convergence}.
However, a great advantage of the ensemble (regression) gradient
is that the blurring takes place implicitly.

Gradient ascent being a local optimisation method,
this implicit smoothing of the objective may well be beneficial --
indeed even a prerequisite --
if the objective possesses many local optima,
as illustrated in \cref{fig:rastrigin},
and thus a reason not to use a very small $\mat{C}_{\vect{u}}$.
Of course, in later iterations, and hopefully near a good optimum,
one could decrease the blurring by decreasing $\mat{C}_{\vect{u}}$.
This could be done with a predefined sequence of scaling values,
by covariance matrix adaptation \citep[CMA,][]{fonseca2015improving},
or by using the gradient with respect to $\mat{C}_{\vect{u}}$ \citep[][]{stordal2016theoretical}.

\subsection{Preconditioning}%
\label{sub:precond}
Some authors prefer to step along a \enquote*{search direction} of
$\ensgrad\! \Loss\, \mat{C}_{\vect{u}}$ rather than the vanilla gradient estimate, $\ensgrad\! \Loss$, itself.
Here, note that $\ensgrad\! \Loss \in \Reals^{1 \times d_{\vect{u}}}$ is a row vector
and that, in contrast with data assimilation,
the covariance matrix, $\mat{C}_{\vect{u}}$, is known,
and might even be small enough to keep in memory.
Alternatively, some prefer $\ensgrad\! \Loss\, \bar{\mat{C}}_{\vect{u}} = \bar{\mat{C}}_{\Loss, \vect{u}}$
where the equality follows by the pseudo-inverse identity \eqref{eq:pinv_identity}.
This preference may be partly due to a reluctance to use pseudo-inversions,
either for computational reasons or for fear of aggravating sampling errors.
However, there is a sense in which sampling errors cancel out for
$\ensgrad\! \Loss\, \mat{C}_{\vect{u}}$,
and not for $\ensgrad\! \Loss\, \bar{\mat{C}}_{\vect{u}}$.
Consider the simple case:
$\Loss(u) = \alpha + \beta u$.
Then $\ensgrad\! \Loss \, C_u = \beta C_u$ exactly, which is not subject to sampling error,
while $\bar{C}_{\Loss, u} = \beta \bar{C}_{u}$ is.
Either way, it is important to mind the potential practical need for
truncating its singular value spectrum or applying Tikhonov regularisation,
detailed in the experiments of \cref{sec:experiments}.

\Chen\ justified the implicit inclusion of $\bar{\mat{C}}_{\vect{u}}$
(and by the same token, $\mat{C}_{\vect{u}}$)
as a regularizing preconditioner for gradient descent.
This claims was further substantiated by the \keyw{natural gradient} perspective
of \citet[][]{stordal2016theoretical}.
Additionally, the design of \Chen\ (and their sometimes \enquote*{double} application) of $\mat{C}_{\vect{u}}$
was used to promote smoothness in the resulting control vector.
However, \citet[][hereafter \Fonseca]{fonseca2017stochastic}
cautioned that the design (e.g. correlation length)
of $\mat{C}_{\vect{u}}$ may not be easy,
and that truncation may complicate the picture.

A potential downside to $\ensgrad\! \Loss\, \mat{C}_{\vect{u}}$
and $\ensgrad\! \Loss\, \bar{\mat{C}}_{\vect{u}}$
is that the direction is no longer estimating the direction of steepest descent,
making the algorithm less transparent.
On the other hand it may be argued that cross-covariances do possess
the essential quality of the gradient, namely that
it (specifically, the signs of its elements)
indicates which coordinate directions are uphill and downhill.
Another clear advantage of $\ensgrad\! \Loss\, \mat{C}_{\vect{u}}$
and $\ensgrad\! \Loss\, \bar{\mat{C}}_{\vect{u}}$,
similar to preconditioning with the actual Hessian,
is that their units contain those of $\vect{u}$, and not its reciprocals,
so that the magnitudes of the individual elements
are likely better suited.

The question of whether and how to precondition
will not be further discussed in this manuscript,
which will focus on $\ensgrad\! \Loss$ alone.
However, all of the equations ending in a trailing $\Ut\pinv$
may be converted to the preconditioned form simply by replacing
$\Ut\pinv$ with $\Ut\tr/(N-1)$,
and most of the discussion is relevant for either choice,
since the question of preconditioning seems independent
of that of dealing with \keyw{robust control}
(also called \keyw{robust optimisation}) problems.

\section{Use in robust control}%
\label{sec:robust_control}
In robust optimisation
the objective function
takes the form of an average, i.e.
\begin{equation}
    \label{eq:Loss}
    \Loss (\vect{u}) \coloneqq \frac{1}{M} \sum_m \loss (\vect{x}_m, \vect{u})\,,
\end{equation}
where the \keyw{conditional objective},
$\loss : \Reals^{d_{\vect{x}} \times d_{\vect{u}}} \rightarrow \Reals$,
is a function also of $\vect{x} \in \Reals^{d_{\vect{x}}}$,
whose uncertainty is quantified by an ensemble
(independent of that of $\vect{u}$) thereof, $\{\vect{x}_m\}$,
where $m=1, \ldots, M$.

Generally, $\loss$ will consist of some utility function
applied to the output of some simulator model for a given $(\vect{x}, \vect{u})$.
For example, preference or net profit, or proximity to a desired, reference, outcome.
The subject of constraints \citep[][e.g.]{oguntola2021ensemble} is not considered,
and the specifics of $\loss$ are not discussed herein.
I.e. we view $\loss$ as a black box,
except to say that EnOpt is mainly tailored for applications in which
evaluations of $\loss$ are costly
because it contains a simulation model
which typically involves computational fluid dynamics.
Indeed, the original application of EnOpt of \Chen\ is petroleum production optimisation,
where $\loss$ is the simulated production of a reservoir,
$\vect{x}$ parametrises its geological uncertainty,
and $\vect{u}$ are well operation control parameters.
Other uses include minimising the wake effect in wind farm layout and operations,
and climate geoengineering studies.

Average objectives like $\Loss$ are also common in statistical inference and machine learning.
In our study, however, $\vect{x}_m$ is not merely a data point,
for which $\loss(\vect{x}_m, \vect{u})$ measures the data mismatch given $\vect{u}$.
Instead, $\vect{x}_m$ is an \emph{input} to the simulator,
i.e. a possible realisation of a model parameter (vector).
This means that evaluating $\loss$
for different values of $\vect{x}$ is equally costly as for different values of $\vect{u}$.
Therefore $M$, i.e. the number of realisations of $\vect{x}$,
tends to be relatively small, usually somewhere between $10$ and $100$,
prompting the utmost parsimony in the use of evaluations for different $\vect{x}$.
For this reason, quasi-Monte-Carlo sampling strategies should be considered
\citep{ramaswamy2020improved}.

\subsection{Analytic gradient}%
\label{sub:analytic_gradient}
In classical steepest (i.e. gradient) descent, iterative, cumulative steps are taken
along the direction of the gradient,
\begin{equation}
    \label{eq:Grad}
    \nabla\! \Loss (\vect{u}) = \frac{1}{M} \sum_m \pd{\loss}{\vect{u}} (\vect{x}_m, \vect{u})
    \,,
\end{equation}
at the current iterate, $\vect{u}$.
The iteration index is not denoted,
as we are only concerned with the accuracy of gradient estimates
for any given iteration.

\subsection{Plain EnOpt gradient}%
\label{sub:basic_enopt_gradient}
The analytic gradient, $\pd{\loss}{\vect{u}}$, is not generally available.
Therefore, in basic EnOpt,
the gradient is estimated by the LLS ensemble gradient \eqref{eq:ensgrad},
yielding
\begin{align}
    \label{eq:ensgrad_Loss}
    \ensgrad\! \Loss
    &=
    \frac{1}{M} \sum_m \loss(\vect{x}_m, \mat{U})\, \Ut\pinv \\
    \label{eq:ensgrad_Loss_lim}
    &\xrightarrow[N]{} \Expect \nabla\! \Loss(\vect{u})
    \,,
\end{align}
where the second line follows from \cref{eq:ensgrad_lim},
and $\loss(\vect{x}_m, \mat{U})$ indicates the row
vector whose $n$-th entry is $\loss(\vect{x}_m, \vect{u}_n)$.

Since the ensemble $\mat{U}$ is random (and resampled at each iteration),
the gradient estimate, $\ensgrad\! \Loss$, is stochastic.
This gives associations to stochastic gradient descent (SGD).
However, in SGD, stochasticity generally arises from a different source:
random sub-sampling of the data, $\{ \vect{x}_m \}$, at each iteration,
to lower costs, to form a \keyw{subgradient} based only on that subset.
Moreover, analytic gradients are often available.
The sub-sampling technique could also be employed with EnOpt,
although the potential benefit seems less, given the relatively small size of $M$.
Nevertheless, in the context of SPSA 
for robust reservoir production optimisation,
\citet[][]{li2013simultaneous,jesmani2020reduced}
did find significant success with sub-sampling.

\section{Using fewer evaluations}%
\label{sec:reducingN}
Computing the plain EnOpt gradient estimate,
$\ensgrad\! \Loss$ of \cref{eq:ensgrad_Loss},
implies $M \times N$ costly evaluations, $\loss(\vect{x}_m, \vect{u}_n)$.

\subsection{The fragile (\enquote*{mean-model}) gradient}%
\label{sub:the_fragile_gradient}
One possibility to effectively reduce $M$ to $1$
is to only evaluate the objective function at $\bar{\vect{x}}$, i.e.
\begin{align}
    \label{eq:grad_fragile1}
    \ensgrad \loss_{|\bar{\vect{x}}}
    &=
    \loss(\bar{\vect{x}}, \mat{U})\, \Ut\pinv \\
    \label{eq:grad_fragile2}
    &\xrightarrow[N]{}
    \Expect \nabla \loss_{|\bar{\vect{x}}}(\vect{u})
    \,,
\end{align}
where the notation $\loss_{|\bar{\vect{x}}}$
indicates the partial function with $\vect{x}$ fixed at $\bar{\vect{x}}$.
The drawback of the estimate $\ensgrad \loss_{|\bar{\vect{x}}}$ is of course that it
does not take into account the uncertainty in $\vect{x}$,
in direct contrast to the express aim of EnOpt
\eqref{eq:ensgrad_Loss},
and may be facetiously termed \keyw{fragile}.

\subsection{The paired gradient}%
\label{sub:the_paired_gradient}
In the original formulation of EnOpt,
inspired by \citet{lorentzen2006new} and \citet{nwaozo2006},
\Chen\ proposed \enquote*{mixing}, or \enquote*{coupling},
or \keyw{pairing} the two ensembles,
meaning that the summation (i.e. averaging) inherent in the matrix product
in \cref{eq:ensgrad_Loss}, which operates over $\vect{u}_n$,
is merged with that over $\vect{x}_m$.
This requires equality of the sample sizes (or that they be divisible),
and only uses $M = N$ evaluations for the gradient estimate,
\begin{equation}
    \label{eq:ensgrad_paired}
    \ensgrad \! \Loss_{\tn{Paired}}
    =
    \loss(\mat{X}, \mat{U})\, \Ut\pinv
    \,,
\end{equation}
where $\mat{X}$ is the ensemble matrix for $\{ \vect{x}_m \}$.
The left-hand side will be explained in \cref{sec:decoupling},
and so should for now be considered merely a label.

In general,
$\Expect_{\vect{u}} \Expect_{\vect{x}} \vect{f}(\vect{x}, \vect{u})
= \Expect_{\vect{x}, \vect{u}} \vect{f}(\vect{x}, \vect{u})$
if $\vect{x}$ and $\vect{u}$ are independent of one another.
Similarly,
\begin{equation}
    \label{eq:nested_and_joint_sample_averages}
    \frac{1}{N} \sum_n \frac{1}{M} \sum_m \vect{f}(\vect{x}_m, \vect{u}_n)
    \approx
    \frac{1}{N} \sum_n \vect{f}(\vect{x}_n, \vect{u}_n)
    \,,
\end{equation}
which becomes an exact equality upon taking the expectation.
\citet[][]{stordal2016theoretical} applied this reasoning to the
plain, unpaired estimate \eqref{eq:ensgrad_Loss}
and the paired estimate \eqref{eq:ensgrad_paired},
both in their preconditioned form, to justify pairing.
They also showed that they share the same, unbiased, expected value,
namely $\Expect \nabla\! \Loss(\vect{u})\, \mat{C}_{\vect{u}}$.
While the pseudo-inverses complicate the picture,
similar reasoning may be applied for the non-preconditioned estimates,
i.e. that the limit of $\ensgrad \! \Loss_{\tn{Paired}}$ is indeed $\Expect \nabla\! \Loss(\vect{u})$.

\subsection{The StoSAG gradient}%
\label{sub:StoSAG}
\citet{fonseca2014robust} proposed 
a \enquote*{modified} EnOpt gradient,
renamed by \Fonseca\ as \keyw{stochastic simplex approximate gradient (StoSAG}).
\begin{equation}
    \label{eq:ensgrad_StoSAG}
    \ensgrad \Loss_{\tn{StoSAG}}
    =
    \bigl[\loss(\mat{X}, \mat{U}) - \loss(\mat{X}, \vect{\mu}) \bigr]\, \Ut\pinv
    \,.
\end{equation}
The second term appears to incur an extra $M$ evaluations of $\loss$.
However, this is not actually an additional cost,
since $\loss(\mat{X}, \vect{\mu})$ and its mean
are normally computed before accepting the previous iterative step.

Since $\loss(\mat{X}, \vect{\mu})$ and $\Ut\tr$ are independent,
and the expectation of $\Ut$ is zero,
the preconditioned form of the estimator \eqref{eq:ensgrad_StoSAG}
is an unbiased estimator just like the paired estimator.
Similarly, the limit of $\ensgrad \Loss_{\tn{StoSAG}}$ is also $\Expect \nabla\! \Loss(\vect{u})$.

Concerning $M < \infty$,
\citet{fonseca2015quantification} and \Fonseca\
have provided considerable empirical evidence of significant improvements in the estimate
\eqref{eq:ensgrad_StoSAG} versus \eqref{eq:ensgrad_paired}.
On the other hand
there exists little theoretical motivation beyond
using a different \enquote{weighting}.
\Fonseca\ did show StoSAG to be superior if $\mat{C}_{\vect{u}}$ is large.
\Cref{sec:fonseca} argues that this conclusion is correct,
but that their analysis is quite flawed.
Our main objective is to show that
the estimate \eqref{eq:ensgrad_StoSAG} hopes to capture
variation due to $\vect{u}$ only,
which is why and how it
improves on the estimate \eqref{eq:ensgrad_paired} of \Chen.

\section{Decoupling the pairing}%
\label{sec:decoupling}
The words \enquote{pairing} and \enquote{coupling} evoke
a deterministic relationship.
Let us explicitly assume such a coupling,
i.e. that $\vect{x}$ is a function of $\vect{u}$,
i.e. $\vect{x} = \vect{\phi}(\vect{u})$.
This assumption immediately yields the effect of
merging the averaging in $\vect{x}$ and $\vect{u}$,
as in \cref{eq:ensgrad_paired},
since $p(\vect{x}|\vect{u}) = \delta(\vect{x} - \vect{\phi}(\vect{u}))$ so that
$\Expect_{\vect{x}, \vect{u}} f(\vect{x}, \vect{u})
= \Expect_{\vect{u}} f(\vect{\phi}(\vect{u}), \vect{u})$.
Indeed, define
\begin{equation}
    \label{eq:Loss_paired}
    \Loss_{\textnormal{Paired}}(\vect{u})
    \coloneqq
    \loss(\vect{\phi}(\vect{u}), \vect{u})
    \,,
\end{equation}
and let $\vect{\phi}$ be such that it
exactly maps each $\vect{u}_n$ to $\vect{x}_n$ for $n = 1,\ldots, M$;
one could for example use a spline interpolant.
Then, since $\Loss_\tn{Paired}(\vect{u}_n) = \loss (\vect{x}_n, \vect{u}_n)$ for all $n$,
applying $\ensgrad$ to $\Loss_\tn{Paired}$
recovers the paired gradient estimate \eqref{eq:ensgrad_paired} of \Chen.

\subsection{The problem with pairing}%
\label{sub:the_problem_with_pairing}
By \cref{eq:ensgrad_lim},
$\ensgrad\! \Loss_\tn{Paired} \xrightarrow[N]{} \Expect \nabla\! \Loss_\tn{Paired}(\vect{u})$.
But
\begin{equation}
    \label{eq:total}
    \nabla\! \Loss_\tn{Paired}
    =
    \pd{\loss}{\vect{u}} + \pd{\loss}{\vect{x}} \frac{\partial \vect{x}}{\partial \vect{u}}
    \,,
\end{equation}
which is the \emph{total} derivative.
By contrast, the plain, unpaired EnOpt gradient \eqref{eq:ensgrad_Loss}
targets \cref{eq:ensgrad_Loss_lim},
i.e. $\frac{1}{M} \sum_m \Expect
\frac{\partial \loss}{\partial \vect{u}} (\vect{x}_m, \vect{u})$,
which contains only the \emph{partial} derivative,
i.e. only the first term of \cref{eq:total},
$\Expect \pd{\loss}{\vect{u}}(\vect{\phi}(\vect{u}), \vect{u})$.

Thus, the introduction of the explicit coupling \emph{function},
which enables the highly articulate tools of analytical calculus,
shows that replacing $\Loss$ with $\Loss_\tn{Paired}$ introduces an error,
represented by the second term of the right-hand side of \cref{eq:total}.
The following subsection corrects for the error,
while the subsequent subsection takes the approach of preventing it.
Both are analytically shown to work in the linear case in \cref{sec:linear_example}.

\subsection{Correcting for the error}%
\label{sub:compensating}
Here we try to compensate
$\Expect \pd{\loss}{\vect{x}}(\vect{\phi}(\vect{u}), \vect{u})
\frac{\partial \vect{x}}{\partial \vect{u}}(\vect{u})$
from $\ensgrad\! \Loss_\tn{Paired}$ and \cref{eq:total}.
Note that this chained derivative is the gradient, evaluated at $\vect{v} = \vect{u}$,
of the partial, composite function $\vect{v} \mapsto \loss_{|\vect{u}} (\vect{\phi} (\vect{v}))$.
Estimating its expected gradient by $\ensgrad$
still leaves the fixing, subscripted $\vect{u}$ unset,
whereas it should also be averaged out.
But since a second sample average would imply $N^2$ evaluations,
we instead make the approximation
\begin{equation}
    \label{eq:loss_at_mean_u}
    \loss_{|\vect{u}} \approx \loss_{|\vect{\mu}}
    \,,
\end{equation}
similar to the approximation of the fragile estimate \eqref{eq:grad_fragile1},
but hopefully less detrimental since is applies to the second(ary), corrective, term.
In summary,
\begin{equation}
    \label{eq:Loss_StoSAG}
    \Loss_\tn{StoSAG}
    :=
    \Loss_\tn{Paired} - \loss_{|\vect{\mu}} \circ \vect{\phi}
    \,,
\end{equation}
and applying $\ensgrad\!$ to $\Loss_\tn{StoSAG}$
recovers the gradient estimate \eqref{eq:ensgrad_StoSAG} of \Fonseca.

\subsection{Preventing the error (decorrelating)}%
\label{sub:decorrelating}

Alternatively, one can eliminate
$\ensgrad \loss_{\vect{\mu}}(\vect{\phi}(\vect{u})) = \loss(\mat{X}, \vect{\mu}) \, \Ut\pinv$,
where we have again included the approximation of \cref{eq:loss_at_mean_u}.
To achieve this, project $\mat{U}$ (row-wise)
onto the space orthogonal to the centred vector $\loss(\mat{X}, \vect{\mu})$,
say $\vect{\psi}$.
The projection is given by
$\mat{U}' = \mat{U} - \mat{U} \vect{\psi} \vect{\psi}\tr / \| \vect{\psi} \|^2_2$,
which we then translate and scale so as to have the same
mean and variances as the original $\mat{U}$.
The gradient is then simply estimated by the paired estimate \eqref{eq:ensgrad_paired},
except with $\mat{U}'$ in place of $\mat{U}$.
A notable drawback of the projection is the potential loss of rank in case $M \le d_{\vect{u}}$.

\section{Generalised StoSAG}%
\begin{figure}[H]
    \centering
    \includegraphics[width=0.99\linewidth]{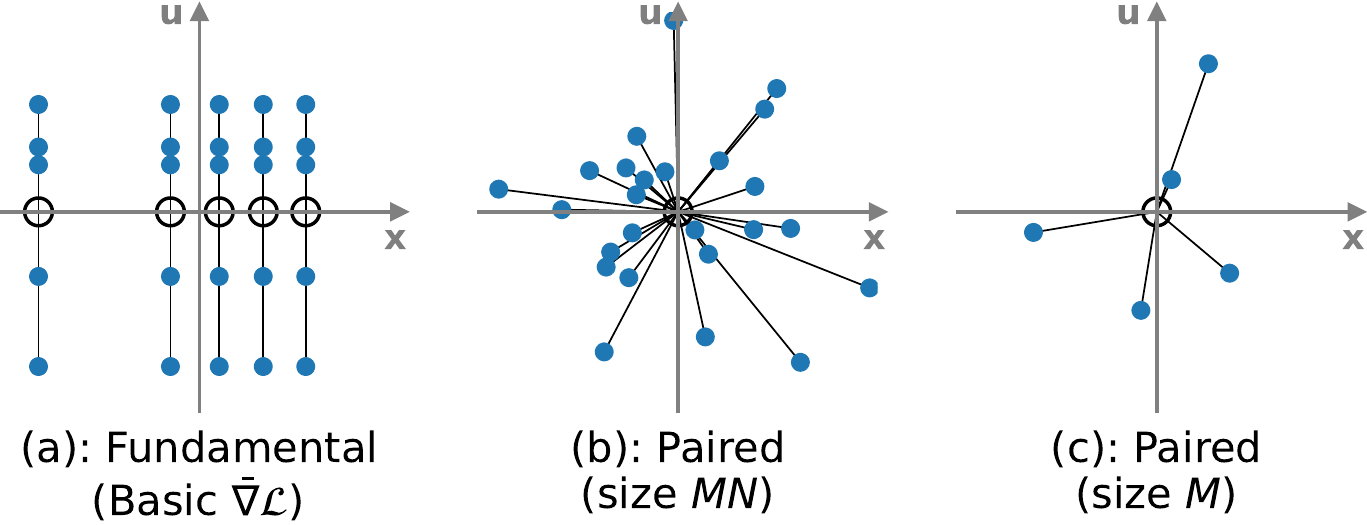}
    \caption{Scatter plot of sample points used in the evaluation of
        $\loss(\vect{x}, \vect{u})$ in various EnOpt gradient estimates.
        The edges and circles indicate the points used to compute increments of $\loss$.
    }
    \label{fig:ux1}
\end{figure}

\label{sec:general_stosag}
As mentioned, the approximation \eqref{eq:nested_and_joint_sample_averages}
may be used to justify the pairing of the gradient estimate \eqref{eq:ensgrad_paired}.
But it can also be used to argue contrarily:
the fact that it is an approximation (for $N < \infty$)
means that paired estimates are not equal to plain, unpaired estimates.
\Cref{fig:ux1} illustrates two issues with pairing:
fewer sample points are used, as intended,
but additionally the sample points are no longer on a grid.
By contrast, the \enquote*{nested} sampling
of the left-hand side of \cref{eq:nested_and_joint_sample_averages},
illustrated by the points in \cref{fig:ux1} (a),
is explicitly independent: the very same ensemble,
$\{ \vect{u}_n \}$, is identically repeated for each $\vect{x}_m$,
which enables showing $\bar{\mat{C}}_{\vect{x},\vect{u}} = \mat{0}$.
Meanwhile, the sample covariance of the points on the right-hand side of
\cref{eq:nested_and_joint_sample_averages}, illustrated in either subfigure (b) or (c),
almost surely contains sampling error,
which make it seem like there is a dependence between $\vect{x}$ and $\vect{u}$.

\begin{figure}[H]
    \centering
    \includegraphics[width=0.99\linewidth]{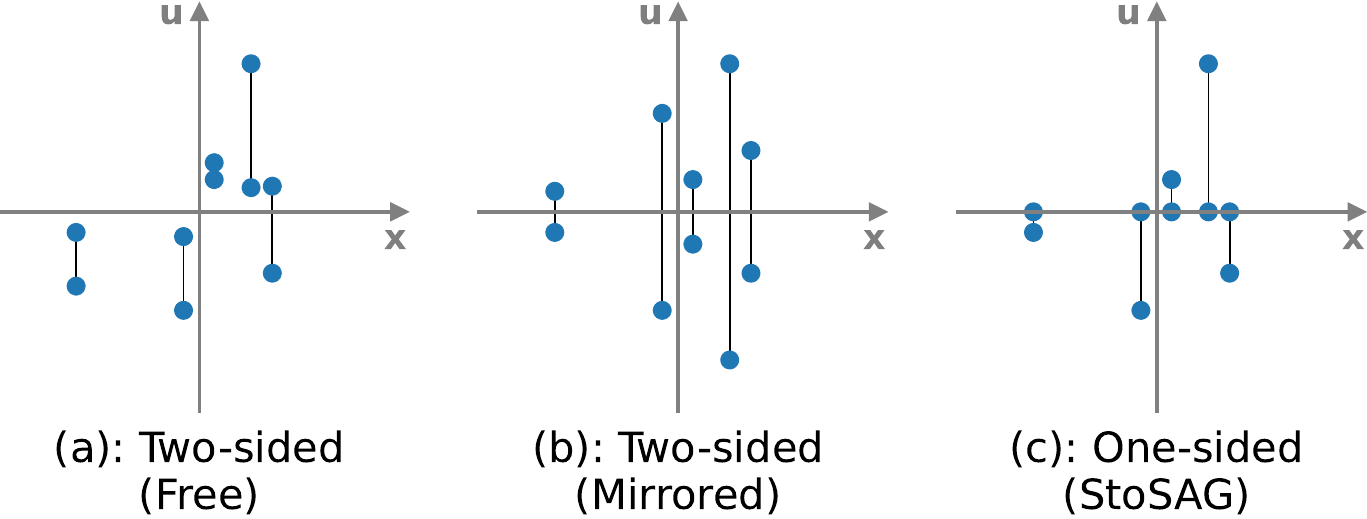}
    \caption{
        Like \cref{fig:ux1}, but for two-sided, distributed ensemble gradients,
        generalizing StoSAG.
    }
    \label{fig:ux2}
\end{figure}

With \cref{fig:ux1} in mind,
the idea of StoSAG gradient is evident.
It is illustrated in \cref{fig:ux2},
which shows three variants that will be developed below,
all of which have $\bar{\mat{C}}_{\vect{x},\vect{u}} = \mat{0}$.

\subsection{Swapping operators}%
\label{sub:swapping_operators}
Note that the plain, unpaired EnOpt gradient \eqref{eq:ensgrad_Loss} may be expressed as
\begin{equation}
    \label{eq:ensgrad_distributed}
    \ensgrad\! \Loss
    =
    \frac{1}{M} \sum_m \ensgrad \loss_{|\vect{x}_m}
    \,,
\end{equation}
i.e., just like the case with the analytic gradient \eqref{eq:Grad},
the ensemble gradient \eqref{eq:ensgrad} distributes over the average in $\vect{x}$.
Thus, we can use a new sample of $\vect{u}$ for each $\vect{x}_m$,
i.e. $\ensgrad \loss_{|\vect{x}_m} = \loss(\vect{x}_m, \mat{U}_m)\, \Ut_m\pinv$.
Then, since $\mat{U}_m$ differs for each $\vect{x}_m$,
we can make its size, $N_m$ (since it can vary with $m$),
small without impoverishing the overall sampling of $\vect{u}$.
In the following experiments, this estimator
will be called \keyw{average LLS} gradient.

It should be noted that equation (12) of \Fonseca\ implicitly takes advantage
of \cref{eq:ensgrad_distributed} [their equation (25)],
but this is not explicitly stated,
nor is their equation (12) motivated except by a \enquote{change of notation}.
Lastly, from our non-preconditioned starting point,
the following additional modification is necessary before arriving at StoSAG.

\subsection{Ratio averaging}%
\label{sub:ratio_averages}
A small $N_m$ necessarily yields significant sampling error;
the pseudo-inversions likely exacerbates it, and causes a significant bias.
Reversing the pseudo-inverse identity \eqref{eq:pinv_identity} yields
\begin{equation}
    \label{eq:ensgrad_m}
    \ensgrad \loss_{|\vect{x}_m} = \bar{\vect{c}}_m\, \bar{\mat{C}}_{\vect{u}, m}\pinv
    \,,
\end{equation}
where $\bar{\mat{C}}_{\vect{u}, m} = \frac{1}{N_m - 1} \Ut_m \Ut_m\tr$
is the $m$-th sample covariance,
and similarly $\bar{\vect{c}}_m = \frac{1}{N_m -1} \loss(\vect{x}_m, \mat{U}_m) \, \Ut_m\tr$.
Now, to alleviate the issues of the pseudo-inversion,
replace $\bar{\mat{C}}_{\vect{u}, m}$ by its average, $\frac{1}{M} \sum_m \bar{\mat{C}}_{\vect{u}, m}$,
so that \cref{eq:ensgrad_distributed} becomes
\begin{align}
    \label{eq:ensgrad_multiplex}
    \ensgrad\! \Loss
    \approx
    \Bigl(\sum_m \bar{\vect{c}}_m\Bigr)
    \Bigl(\sum_m \bar{\mat{C}}_{\vect{u},m}\Bigr)\pinv
    \,.
\end{align}
In other words, the average of the \enquote*{ratio} \eqref{eq:ensgrad_distributed}
is replaced by the ratio of averages \eqref{eq:ensgrad_multiplex}.
This forms the generalised StoSAG ensemble gradient,
for which different sampling strategies (all using $N_m=2$) are illustrated in \cref{fig:ux2}.

Since covariance estimates are themselves averages,
it is tempting to also replace
$\frac{1}{M} \sum_m \bar{\vect{c}}_m$,
by the plain covariance estimate computed by pooling all of the sub-samples,
and do likewise for the covariances.
However, doing so would merge the averaging in $\vect{u}$ with that in $\vect{x}$,
and indeed would just recover the paired gradient estimate \eqref{eq:ensgrad_paired} of \Chen.
Alternatively, one can do so for the covariances alone,
resulting in an estimator we call \keyw{hybrid}.

\subsection{Using only two points}%
\label{sub:using_n_2_}
A particularly interesting case is $N_m=2$,
in which case we
denote the \emph{two} columns of $\mat{U}_m$
by $\vect{v}_m$ and $\vect{w}_m$.
Then, since $\Ut_m$ is centred,
its two columns,
$\tilde{\vect{v}}_m$ and $\tilde{\vect{w}}_m$,
are symmetric about zero, i.e.
$\tilde{\vect{v}}_m = - \tilde{\vect{w}}_m = ( \vect{v}_m - \vect{w}_m )/2$,
so that
$\bar{\mat{C}}_{\vect{u}, m} = \Ut_m \Ut_m\tr
=
\tilde{\vect{v}}_m \tilde{\vect{v}}_m\tr +
\tilde{\vect{w}}_m \tilde{\vect{w}}_m\tr
= 2 \tilde{\vect{v}}_m \tilde{\vect{v}}_m\tr
$.
Hence
\begin{equation}
    \label{eq:mean_cov2}
    \sum_m \bar{\mat{C}}_{\vect{u},m}
    =
    \frac{1}{2} (\mat{V} - \mat{W}) (\mat{V} - \mat{W})\tr
    \,,
\end{equation}
where the $m$-th columns of $\mat{V}, \mat{W} \in \Reals^{d_{\vect{u}} \times M}$
are $\vect{v}_m$ and $\vect{w}_m$.
Similarly,
\begin{equation}
    \label{eq:crosscov_m}
    \bar{\vect{c}}_m
    =
    \frac{1}{2} \bigl[
        \loss(\vect{x}_m, \vect{v}_m) -
        \loss(\vect{x}_m, \vect{w}_m)
    \bigr] \, ( \vect{v}_m - \vect{w}_m )\tr
    \,,
\end{equation}
whose sum in $m$ can be written as a matrix product
similar to \cref{eq:mean_cov2}.
Thus, \cref{eq:ensgrad_multiplex} becomes
\begin{equation}
    \label{eq:ensgrad_2s}
    \ensgrad_{\!\tn{Two-sided}}\, \Loss
    \coloneqq \Bigl[ \loss(\mat{X}, \mat{V}) - \loss(\mat{X}, \mat{W}) \Bigr]
    (\mat{V} - \mat{W})\pinv
    \,,
\end{equation}
where we have again used the pseudo-inverse identity \eqref{eq:pinv_identity}.
In summary, in the case of $N_m = 2$,
we can write \cref{eq:ensgrad_multiplex}
in a manner reminiscent of the ensemble LLS gradient \eqref{eq:ensgrad}
as well as two-sided finite differences.

Also note that in this $N_m = 2$ case,
the pseudo-inverses, $\{\Ut_m\pinv\}$, in the average LLS gradient
\eqref{eq:ensgrad_distributed} are analytically available
and so can be efficiently computed.

\subsection{Using mirror points}%
\label{sub:mirrored}
An alternative to sampling $\{ \vect{w}_m \}$
is to simply set $\vect{w}_m$ as the reflection of $\vect{v}_m$ about the mean,
i.e. $\vect{w}_m - \vect{\mu} = \vect{\mu} - \vect{v}_m$.
This may be seen as a lightweight sampling strategy,
in the vein of quasi-random or quasi-Monte-Carlo methods.
Of course, covariance estimators should include a scaling of $1/2$
when using mirrored samples,
but this effect cancels in \cref{eq:ensgrad_multiplex}.

With \keyw{mirroring}, both $\vect{v}_m$ and $\vect{w}_m$
may be written in terms of a single draw, say $\vect{u}_m$ (again),
i.e. $\mat{V} = \vect{\mu} + \Ut$
and $\mat{W} = \vect{\mu} - \Ut$,
where the $m$-th column of $\Ut$ is $\vect{u}_m - \vect{\mu}$,
and the vector-matrix sum is understood appropriately.
Thus, $\mat{V} - \mat{W} = 2 \Ut$, and
\begin{equation}
    \label{eq:ensgrad_mirrored}
    \ensgrad_{\!\tn{Two-sided}}\, \Loss
    = \frac{1}{2} \Bigl[
        \loss(\mat{X}, \vect{\mu} {+} \Ut) - \loss(\mat{X}, \vect{\mu} {-} \Ut)
    \Bigr]
    \, \Ut\pinv
    \,.
\end{equation}

\subsection{Mirroring values}%
\label{sub:one_sided}
Finally, in view of the above symmetry,
it particularly tempting to make the approximation
\begin{equation}
    \label{eq:loss_mirrored}
    \loss(\mat{X}, \vect{\mu})
    \approx
    \frac{1}{2} \bigl[\loss(\mat{X}, \vect{\mu} {-} \Ut)
    + \loss(\mat{X}, \vect{\mu} {+} \Ut)\bigr]
    \,.
\end{equation}
Solving \cref{eq:loss_mirrored} for
$\loss(\mat{X}, \vect{\mu} {-} \Ut)$,
yields the reflection of $\loss(\mat{X}, \vect{\mu} {+} \Ut)$
about $\loss(\mat{X}, \vect{\mu})$,
and its substitution into \cref{eq:ensgrad_mirrored}
yields a \keyw{one-sided} variation,
which is nothing but
$\ensgrad \! \Loss_{\tn{StoSAG}}$ \eqref{eq:ensgrad_StoSAG}.

The averaging in \cref{eq:loss_mirrored}
recalls the method of \citet{jeong2020efficient} of
\enquote{approximating the objective function values of unperturbed control
variables using the values of perturbed ones}.
However, mirroring avoids the complication of estimating
clusters of $\{ \vect{u}_n \}$ for subset averaging.
Also recall that \cref{eq:ensgrad_multiplex} may be applied
with any $N_m \ge 2$, even differing with $m$ \citep[][]{fonseca2015quantification},
so that larger subsets may be used.

The main advantage of $\ensgrad\! \Loss_\tn{StoSAG}$
is that it uses $\loss(\mat{X}, \vect{\mu})$,
which is available from the \emph{previous} step validation,
whereas the true two-sided variants
require evaluation an additional $N$ points.
Of course, the true two-sided variants can approximate $\loss(\mat{X}, \vect{\mu})$
using $\loss(\mat{X}, \mat{V})$ and $\loss(\mat{X}, \mat{W})$
and, as shown in \cref{fig:ux2},
they yield the benefit of using twice as many values of $\vect{u}$
in computing the gradient.

\begin{figure}[H]
  \centering
  \includegraphics[width=\linewidth]{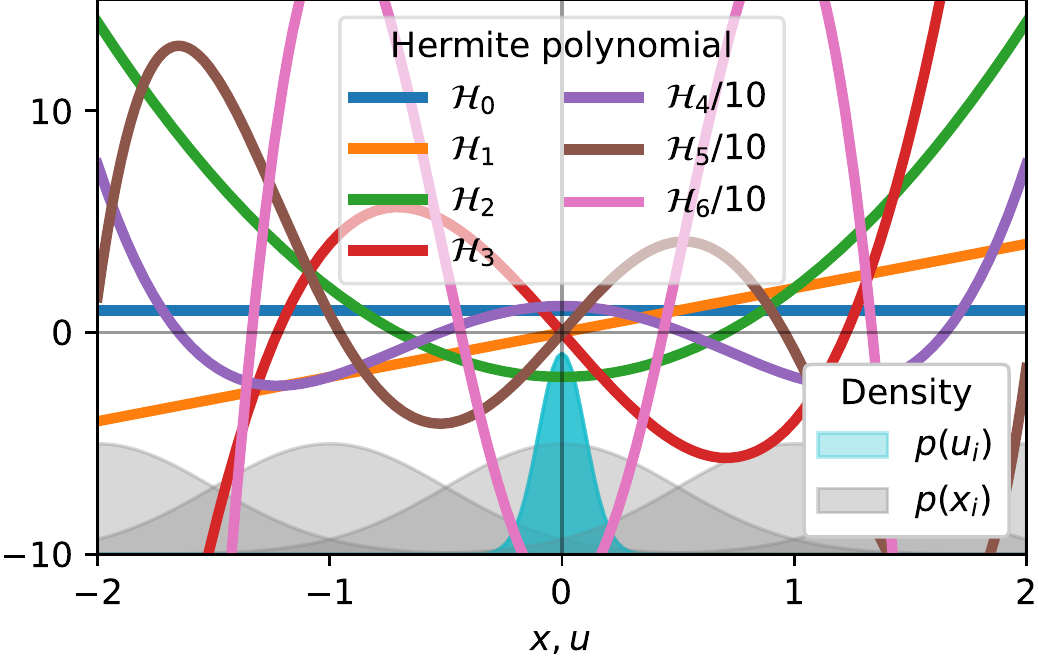}
  \caption{Elements of the experiment.}
  \label{fig:elements}
\end{figure}

\section{Experimental comparison}%
\label{sec:experiments}

\citet[][]{fonseca2015quantification} studied the accuracy of the (preconditioned) StoSAG
gradient estimator on the Rosenbrock function,
comparing it with the paired ensemble gradient.
This section presents a related quantitative comparison,
but including all of the various gradient estimators presented in this paper.

\subsection{Description}%
\label{sub:experiment_description}
The experiment consists of computing the various ensemble gradients
of $\Loss (\vect{u}) = \sum_i \mathcal{H}(u_i + x_i)$,
where $\mathcal{H}$ is a given Hermite polynomial,
$i = 1, \ldots, 5$ is the dimension index,
$\vect{u} \sim \NormDist(\vect{0}, \mat{I}/100)$,
and
$\vect{x} \sim \NormDist( [-2, \ldots, 2]\tr , \mat{I}/4)$,
as illustrated in \cref{fig:elements}.

The error of a gradient estimate is defined as the difference with
the estimation target (\enquote{truth}), $\Expect \nabla\! \Loss(\vect{u})$.
Each experiment is repeated $10^6$ times using different random seeds.
The root-mean-square (RMS) error is then computed over the seeds.
The bias magnitude is obtained by averaging the errors over the seed,
and taking its absolute value.
Subsequently, both bias and RMS error are averaged over the dimensions, $i$.

The pseudo-inverses are all computed using Tikhonov regularisation,
which we found to be more robust than truncation.
Specifically, a pseudo-inverse is computed via singular value decomposition,
and instead of using the reciprocals of the singular values, $1/s_i$,
the diagonal matrix consists of the elements $s_i / \bigl(s_i^2 + (\lambda s_1)^2\bigr)$,
where $\lambda$ is the (tunable) regularisation parameter, and $s_1$ the largest singular value.
A range of $\lambda$ is tried out,
and only the best result of each method is retained.

For generality and sensitivity analysis,
this is all repeated for a range of different ensemble sizes
and for each Hermite polynomial from $0$ to $6$.

The experiments are coded in Python,
using NumPy \citep[][]{harris2020array},
SciPy \citep[][]{2020SciPy-NMeth},
pandas
\citep[][]{mckinney-proc-scipy-2010},
and Xarray
\citep[][]{hoyer2017xarray},
and can be found in the supplementary materials of the paper.

\begin{figure*}[ht]
  \centering
  \includegraphics[width=\textwidth]{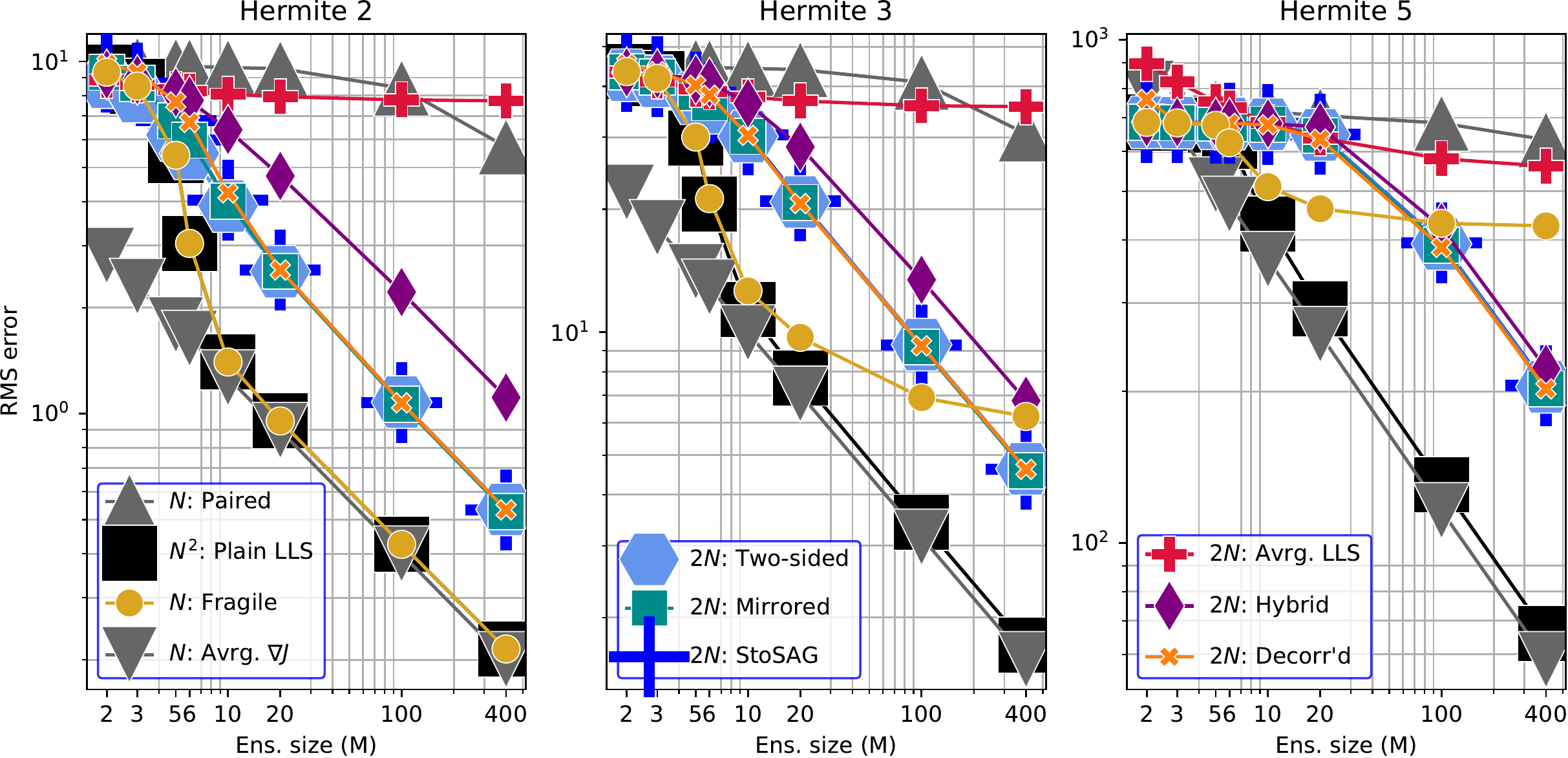}
  \caption{Root-mean-square error as a function of ensemble size
      for different ensemble gradient estimators.
      The panels show experiments using different Hermite polynomials,
      but all legends apply to all panels,
      with the number of evaluations of $\loss$
      (i.e. $N$, $2 N$, $N^2$) shown in front of their labels.
  }
  \label{fig:rmse}
\end{figure*}

\subsection{Results}%
\label{sub:experiment_results}
The RMS errors are shown in \cref{fig:rmse}.
As expected, the errors all increase from left to right across panels (notice the shifting y-axis),
i.e. with increasing Hermite order (a proxy for nonlinearity),
while they decrease inside each panel,
i.e. for increasing ensemble size (on the x-axis).
The \keyw{average $\nabla \Loss$}
(evaluated and averaged over the same ensemble as the other estimators)
is uniformly the best.
It is generally closely followed by the \keyw{plain LLS} gradient estimate,
$\ensgrad \Loss$ of \cref{eq:ensgrad_Loss},
which uses $N^2$ evaluations,
while the paired ensemble gradient,
$\ensgrad \Loss_{\!\tn{paired}}$ of \cref{eq:ensgrad_paired},
is nearly uniformly the worst.

The other methods generally score RMS errors somewhere in between these \enquote*{baselines},
although simply taking the average of LLS ensemble gradients (red)
performs very poorly.
It can be seen that all of the StoSAG variants (blue),
as well as the decorrelated approach (orange), perform almost equally well
across the board.
They appear to attain asymptotic convergence rates almost as soon as $N > 5 = d_{\vect{u}}$,
in which setting the ensemble gradient is computed with a full-rank ensemble.
The scores of the fragile (\enquote*{mean-model}) ensemble gradient (gold) undergoes
a more precipitous drop around $N = 6$,
but flattens out for large ensemble sizes, since it is not a consistent
estimator of $\Expect \nabla\! \Loss(\vect{u})$.
However, for moderate ensemble sizes the StoSAG ensemble gradients
are soundly outperformed by the fragile ensemble gradient,
which means that, in terms of RMS error,
their additional sampling error is more significant
than the bias in the fragile gradient.
The case of the 4th and 6th order Hermite polynomials (not shown)
are similar to that of the 3rd and 5th, the general picture being that
the higher the Hermite order is, the lower (but still high) $N$ is
where the StoSAG and fragile score lines intersect.

Interestingly, for Hermite polynomial of degree $2$ (leftmost panel)
the fragile gradient achieves equally good scores as the plain LLS ensemble gradient
for all ensemble sizes, which is explained by the fact that in this case $\nabla \Loss$ is linear.
In the case of the 1st order Hermite polynomial (not shown)
the gradient is constant, and all of the methods are equally good,
except for the paired ensemble gradient.
In the case of the 0th order (constant, not shown) they
all essentially achieve zero error.

\begin{figure}[H]
  \centering
  \includegraphics[width=\linewidth]{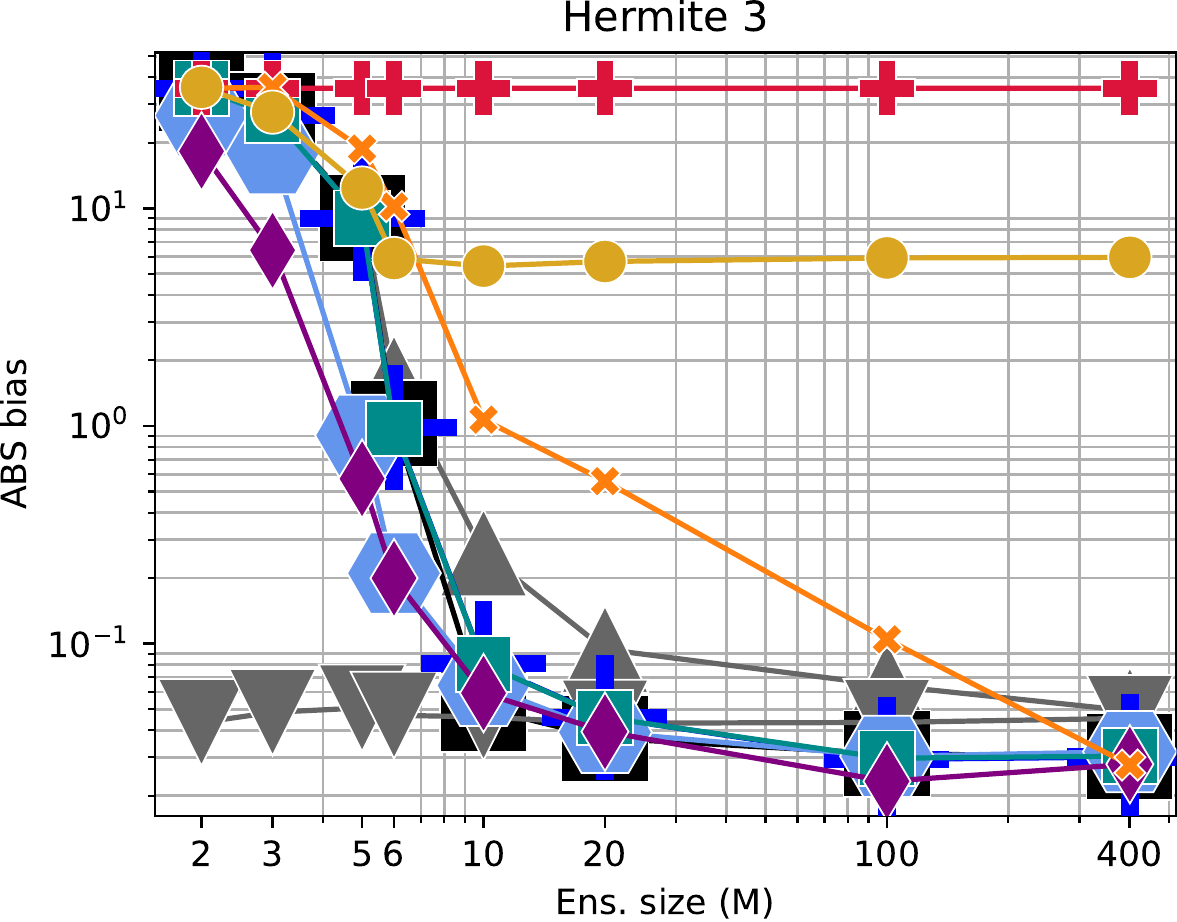}
  \caption{Like \cref{fig:rmse} but for absolute bias
      (note: bias-minimising Tikhonov regularisation values
      generally differ from those used in \cref{fig:rmse}).
  }
  \label{fig:bias}
\end{figure}

The bias is shown in \cref{fig:bias}.
Here, the paired ensemble gradient does not demonstrate a markedly worse bias than StoSAG.
Interestingly, the performance of the hybrid method (purple),
where the covariance average is replaced by the total covariance,
is actually slightly lower than for StoSAG,
whereas in terms of RMS error it was slightly worse.
The decorrelation method exhibit a somewhat higher bias than StoSAG.
Most striking is that the fragile gradient has a very significant bias,
explaining almost all of its RMS error.

\section{Summary}%
\label{sec:conclusion}
As detailed in the abstract and introduction,
the aim of this paper has been to
shed more light on the workings of the ensemble gradient in EnOpt;
in particular its specialisation to robust control problems,
and the importance of decoupling the partial derivatives,
as is done by the StoSAG variant(s).
The numerical experiments indicate that a very large ensemble size
is needed for these robust ensemble gradient estimators
to achieve lower errors than the fragile (\enquote*{mean-model}) estimator,
albeit not in terms of bias.

\subsection*{Acknowledgements}
This work has been funded by the Research Council of Norway
under the CoRea project (grant number 301396),
with some additional funding from Remedy (grant number 336240)
We also thank
Kjersti Solberg Eikrem for pointing out the possibility of centring by $\vect{f}(\vect{\mu})$
in covariance estimation.

\appendix
\section{Centring}%
\label{sec:centring}
We resume the discussion directly below \cref{eq:cov_estim_true_mean}.
The reason its estimate \emph{usually} has a larger variance than \cref{eq:ens_cov}
is of course that \cref{eq:cov_estim_true_mean}
does not centre $\vect{f}(\vect{u}_n)$,
by its sample mean, much less the true mean.
Meanwhile, using the sample mean to centre $\vect{f}(\vect{u}_n)$
in \cref{eq:cov_estim_true_mean} simply recovers $\bar{\mat{C}}_{\vect{f},\vect{u}}$,
except scaled by $\frac{N-1}{N}$ which results in a bias.

Consider, for example, the simple case $f(u) = \beta$,
for some constant $\beta$:
\cref{eq:ens_cov} yields exactly $\bar{C}_{f, u} = 0$, which has zero variance,
while \cref{eq:cov_estim_true_mean} produces $\beta (\bar{u} - \mu)$,
which has variance $\beta C_{f,u}/N$.

In the linear case: $f(u) = \alpha u + \beta$,
it can be shown that
\cref{eq:ens_cov} yields $\bar{C}_{f, u} = \frac{\alpha}{N-1} \sum_n (u_n-\bar{u})^2 $,
while \cref{eq:cov_estim_true_mean} produces
$\frac{\alpha}{N} \sum_n (u_n-\mu)^2 + (\alpha \mu - \beta) (\bar{u} - \mu)$.

Note that, in both cases, 
\cref{eq:cov_estim_true_mean}
produced the additional term $\eta\, (\bar{u} - \mu)$,
where $\eta = \Expect f(u)$,
is the true (but unknown) mean of $f$,
i.e.
$\eta = \beta$ and $\eta = (\alpha \mu - \beta)$, respectively.
Thus, unless $\eta$ happens to be sufficiently near zero,
$\bar{\mat{C}}_{\vect{f}, \vect{u}}$ of \cref{eq:ens_cov}
will generally be less noisy lower than
the estimate of \cref{eq:cov_estim_true_mean}.

\section{Comment on \texorpdfstring{\Fonseca}{Fonseca'2017}}%
\label{sec:fonseca}
As mentioned in \cref{sub:StoSAG},
\Fonseca\ developed a theoretical error analysis
of the preconditioned (i.e. cross-covariance) form of
$\ensgrad \! \Loss_{\tn{paired}}$ of \Chen\,
and $\ensgrad \! \Loss_\tn{StoSAG}$ of \citet{fonseca2014robust}.
However, there are a few issues, detailed below.

Firstly, their preoccupation with approximations (8) and (9) of \Fonseca,
which \Chen\ used in a Taylor development of the cross-covariance,
is misplaced.
\Fonseca\ appears to make the odd claim that the approximation (9)
of the average of $\loss$ by its evaluation \emph{at the average}
requires $\loss(\vect{x}_i, \vect{u}) \approx \loss(\vect{x}_j, \vect{u})$ for any $i,j$,
but it actually only requires a linear assumption.
Moreover, \Chen\ only used it for analysis,
not for constructing the method,
making its importance less \enquote{imperative}.

Actually, the main issue in the derivation of equation (14) of \Chen\ 
is an unstated neglect of the variation in the uncertain parameters
(following their substitution of the Taylor expansion).
However, \Fonseca\ appear to make a very similar omission:
the second term on the right-hand side of their equation (28)
is overlooked in the subsequent analysis.
But this is the leading-order increment due to change in $\vect{x}$,
and the fact that (all) increments due to variation in $\vect{x}$ are zero for StoSAG
is the crucial difference.

Lastly, \Fonseca\ derive the magnitudes of the errors
of the cross covariance estimator, namely $\mathcal{O}{(\max_n \| \vect{u}_n - \vect{\mu} \|^k)}$
with $k=3$ for the StoSAG gradient and only $k=1$ for the paired one.
However, they do not comment on the fact that the cross-covariance
itself scales with $\| \vect{u} - \vect{\mu} \|^2$.
The root problem, moreover, is that a focus on individual (single-member) errors
cannot hope to show that errors cancel by averaging, as pointed out in \cref{sub:stein}.

\section{StoSAG as variance reduction}%
\label{sec:variance_reduction}
The principle of variance reduction can be stated as follows.
Suppose we possess an estimate $a$,
and some other random variable, $b$, whose mean is zero.
As is well known,
\begin{equation}
    \label{eq:var_reduction}
    \Var(a - b) = \Var(a) + \Var(b) - 2 \Cov(a, b)
    \,.
\end{equation}
Thus, while the mean is left unchanged,
the variance of the \enquote*{corrected} estimate $a - b$
might be lower than that of $a$
if the covariance is sufficiently large.
Indeed, rearranging \cref{eq:var_reduction}
it can be seen that the improvement, $\Var(a) - \Var(a - b)$,
is positive if the linear correlation of $a$ and $b$ satisfies
$\rho \ge \frac{1}{2 r}$,
where $r \coloneqq \sqrt{\Var(b) / \Var(a)}$ is the ratio of the variances.

A variance reduction perspective was also considered by \citet[][]{stordal2016theoretical},
based on \citet[][]{sun2009efficient}.
However, they considered the \emph{design} of a zero-mean $b$
so as to yield reduced variances.
This is not generally possible, unless $\Loss$ has a simple expression,
while empirical estimates were not found successful.
By contrast, our purpose here is to use \cref{eq:var_reduction}
to \emph{explain} the success of the StoSAG gradient.

Let
$\vect{a} = \ensgrad\! \Loss_\tn{Paired} = \loss(\mat{X}, \mat{U}) \Ut\pinv$
and
$\vect{a} - \vect{b} = \ensgrad\! \Loss_\tn{StoSAG}$,
where the correction is
$\vect{b} = \loss(\mat{X}, \vect{\mu})\, \Ut\pinv$.
Note that indeed $\Expect \vect{b} = \vect{0}$,
so that the subtraction of $\vect{b}$
does not alter the expected value.
At the same time, $\vect{b}$ may well be expected to co-vary
with $\vect{a}$, since they both contain $\mat{X}$.
Therefore its subtraction from $\vect{a}$ should yield a reduction in variance.

Let us investigate this in the framework of \cref{eq:var_reduction}.
Since $\vect{a}$ and $\vect{b}$ are now vectors,
the variances in \cref{eq:var_reduction} get replaced by covariance matrices,
while $2 \Cov(a, b)$ is replaced by $\Cov(\vect{a}, \vect{b})$ plus its transpose.
Rephrasing the previous paraphrase,
it seems reasonable (but we furnish no guarantees),
that for most $\loss$, the correlation of
$\loss(\mat{X}, \mat{U})$ and $\loss(\mat{X}, \vect{\mu})$
will be significantly positive
($\loss(\mat{X}, \vect{\mu})$ having correlation $1$ with itself,
it difficult to see how replacing $\vect{\mu}$ with $\mat{U}$
would drastically lower it).
Similarly, $\vect{a}$ and $\vect{b}$ should have significantly positive correlation.
Thus, by \cref{eq:var_reduction},
$\ensgrad\! \Loss_\tn{StoSAG}$
should have a lower variance than $\ensgrad\! \Loss_\tn{Paired}$.

By \cref{eq:var_reduction},
the relative improvement $[\Var(a) - \Var(a - b)]/\Var(a)$
is given by $r(2 \rho - r)$.
Thus, if the variance of $\vect{x}$ is small,
then the variance of $\loss(\mat{X}, \vect{\mu})$ is small,
so that the variance of $\vect{b} = \loss(\mat{X}, \vect{\mu})\, \Ut\pinv$
is often small, compared to that of $\vect{a}$.
Hence $r$ is small and so is the relative improvement.

As an example, consider the bilinear conditional objective,
$\loss(\vect{x}, \vect{u}) = \vect{x} + \vect{u}$,
generalisable to higher dimensions by summing over the dimensions.
In this linear case,
\begin{equation}
    \ensgrad\! \Loss_\tn{StoSAG}
    = (\mat{X} + \mat{U})\, \Ut\pinv - (\mat{X} + \vect{\mu}\ones\tr)\, \Ut\pinv
    \,,
\end{equation}
where the first and penultimate terms cancel,
and $\ones\tr \Ut\pinv = \mat{0}$,
since $\Ut\pinv = \Ut\tr (\Ut \Ut)\pinv$ is orthogonal to the vector of ones, $\ones$.
This leaves only $\Ut \Ut\pinv$ which,
if $N > d_{\vect{u}}$, is the centred identity matrix,
which has zero variance.
The same conclusion is provided by \cref{eq:var_reduction}.
In this case, the variances are equal, i.e. $r=1$, and $\rho = 1$,
so that the relative improvement is $1$.
If $N \le d_{\vect{u}}$, then $\Ut \Ut\pinv$ is a projection matrix
(whose eigenvalues are all $0$ or $1$) whose subspace,
reflecting the fact that not the entire space is sampled, is random.
In either case, $\ensgrad\! \Loss_\tn{StoSAG}$ yields
a strong improvement on $\ensgrad\! \Loss_\tn{Paired}$,
which must also contend with the spurious noise of $\mat{X} \Ut\pinv$.

\section{Linear example}%
\label{sec:linear_example}
Consider the simple the bilinear conditional objective,
\begin{equation}
    \label{eq:Loss_linear}
    \loss(\vect{x}, \vect{u})
    = \sum_i (\mat{A} \vect{x} + \mat{B} \vect{u})_i
    = \ones\tr (\mat{A} \vect{x} + \mat{B} \vect{u})
    \,,
\end{equation}
whose gradient is the constant $\pd{\loss}{\vect{u}} = \ones\tr \mat{B}$.
Thus, the target, i.e. the average gradient,
$\Expect \nabla\! \Loss(\vect{u})$, of \cref{eq:ensgrad_Loss_lim},
is also $\ones\tr \mat{B}$.

The paired ensemble gradient is
\begin{equation}
    \label{eq:ensgrad_loss_linear}
    \ensgrad\! \Loss_\tn{Paired}
    = \ones\tr \mat{A} \mat{X} \Ut\pinv + \ones\tr \mat{B} (\mat{U} \Ut\pinv)
    \,.
\end{equation}
But, providing $N - 1 \ge d_{\vect{u}}$,
the rank of $\Ut$ is $d_{\vect{u}}$ and hence $\Ut \Ut\pinv = \mat{I}$.
Thus, since the rows of $\mat{U} - \Ut$ are orthogonal to those of $\Ut$,
it follows by \cref{eq:pinv_identity} that also $\mat{U} \Ut\pinv = \mat{I}$.
In summary, the error in $\ensgrad\! \Loss_\tn{Paired}$
is simply $\ones\tr \mat{A} \mat{X} \Ut\pinv$,
which can be seen to be due to the correlation between $\mat{X}$ and $\mat{U}$.

Meanwhile, as shown by \eqref{eq:Loss_StoSAG},
$\ensgrad_{\!\tn{StoSAG}} \Loss$ subtracts
$\ensgrad (\loss_{|\vect{\mu}} \circ \vect{\phi})
= \ones\tr \mat{A} \mat{X} \Ut\pinv
+ \ones\tr \mat{B} \vect{\mu} \ones\tr \Ut\pinv$.
But since $\Ut$ is centred,
\cref{eq:pinv_identity} yields $\ones\tr \Ut\pinv = \vect{0}$.
Thus, StoSAG obtains the exact gradient.

Similarly, centring $\loss(\mat{X}, \vect{\mu})$ yields
$\vect{\psi} = \ones\tr \mat{A} \tilde{\mat{X}}$.
Thus, the null space of the decorrelated $\mat{U}'$
contains $\ones\tr \mat{A} \tilde{\mat{X}}$,
and substituting $\mat{U}'$ for $\mat{U}$ in \cref{eq:ensgrad_loss_linear}
eliminates the error.

\bibliographystyle{myabbrvnat}
\bibliography{references}
\end{multicols}
\end{document}